\documentclass[11pt]{article}

\usepackage{amsmath} 
\usepackage{amssymb}  
\usepackage{amsthm}

\usepackage{amsfonts}
\usepackage{amscd}
\usepackage{pifont}
\usepackage{mathrsfs}
\usepackage{mathtools}

\usepackage{amsthm}
\usepackage{color}
\usepackage{graphicx}
\usepackage{graphics}
\usepackage{xcolor}
\usepackage{comment}

\newtheorem{theorem}{Theorem}
\newtheorem{lemma}[theorem]{Lemma}
\newtheorem{proposition}[theorem]{Proposition}

\theoremstyle{remark}
\newtheorem{remark}{Remark}

\theoremstyle{definition}
\newtheorem{definition}{Definition}

\newcommand{\N}{\mathbb{N}}
\newcommand{\R}{\mathbb{R}}

\newcommand{\spazio}{\hspace{1mm}}

\newcommand{\diff}{\mathrm{Diff}}

\newcommand{\nota}[1]{#1}

\newcommand{\Diff}{\mathrm{Diff}}
  
\def\vrexp{\stackrel{\longrightarrow}{\exp}}
\newcommand{\rexp}[3]{\,  \vrexp \int_0^{#1} #2 \, d #3}

\author{
Marco Caponigro\thanks{Dipartimento di Matematica, Universit\`a di Roma ``Tor Vergata'', Via della Ricerca Scientifica 1, 00133 Roma, Italy \texttt{caponigro@mat.uniroma2.it}}
\thanks{Corresponding author}
\and
Arianna Vicari \thanks{Dipartimento di Matematica, Universit\`a di Roma ``Tor Vergata'', Via della Ricerca Scientifica 1, 00133 Roma, Italy \texttt{vicari@mat.uniroma2.it}}
}

\title{Transport maps as flows of control-affine systems}

\begin{document}

\date{}
\maketitle

\begin{abstract}
We consider the problem of transporting \nota{one probability measure into another through} the flow of a given driftless control-affine system.
Under suitable regularity conditions, the controllability of the system by means of open-loop controls is a sufficient condition for the existence of time-varying
feedback controls such that the $1$-time flow of the system is the optimal transport map for the quadratic cost.
\end{abstract}

\section{Introduction}
In this paper we study the transport problem under dynamical constraints, more precisely given two probability measures $\mu$ and $\nu$ on an $n$-dimensional oriented Riemannian manifold $(M,g)$, and a driftless control-affine system on $M$
\begin{equation}\label{eq:driftless}
\dot x = \sum_{i=1}^{m} u_{i}f_{i}(x),
\end{equation}
\nota{our goal} is to couple $\mu$ and $\nu$ with the flow, at fixed time, of system~\eqref{eq:driftless}.
\nota{Moreover, in the Euclidean case, we realize the optimal transport map for the quadratic cost $c(x,y)=|x-y|^2$ between $\mu$ and $\nu$, as the flow, at fixed time, of system~\eqref{eq:driftless} associated with suitable time-varying feedback controls $u_1, \dots, u_m$.}

In particular, we focus on the smooth, i.e. $C^\infty$ case, when $\mu$ and $\nu$ are absolutely continuous measures with respect to the Riemannian volume form $dV_g$ \nota{(the Lebesgue measure in the Euclidean setting)} and their densities are smooth and bounded away from zero.

\nota{Our contribution is twofold. On the one hand, when $M$ is a compact Riemmanian manifold, 
we prove that the controllability of system~\eqref{eq:driftless} by means of open-loop controls is a sufficient condition to transport $\mu$ to $\nu$ via the flow of the systems associated with a time-varying feedback.
On the other hand, when $M$ is an open subset of $\R^n$, under the same controllability assumptions on the system, we prove that there exist time-varying feedback controls realizing the optimal transport map between $\mu$ and $\nu$ as the $1$-time flow of the system.
These two results, although independent, share the same techniques in their proofs and they are naturally related since the problem of realizing the optimal transport map can be viewed as the \emph{optimal} version of the coupling problem. }

If $M$ is compact, $m = n$, and the vector fields $f_1, \dots, f_n$ are linearly independent on $M$ (i.e. any direction is admissible), \nota{the problem of coupling $\mu$ and $\nu$} is classical and well known. Indeed Moser's Theorem~\cite{moser} guarantees the existence of a smooth \nota{and isotopic to the identity} diffeomorphism $T$, coupling $\mu$ and $\nu$ in the sense that
$$
\nu(A) = \mu(T^{-1}(A)),
$$
for any $A$ Borel subset of $M$.
From Thurston's result~\cite{thurston1974foliations} on the simplicity of the group of diffeomorphisms, every diffeomorphism isotopic to the identity can be decomposed as a finite product of flows of the vector fields $f_1, \dots, f_n$ rescaled by suitable smooth functions, implying the existence of functions $(t,x) \mapsto u_i(t,x)$, $i =1, \dots, n$, such that $T$ is the flow at time $1$ of system \eqref{eq:driftless} associated with $u_1,\dots,u_n$.

In \nota{this work} we consider the \nota{problem of coupling $\mu$ and $\nu$} in the case in which $m<n$, under the sole assumption that system~\eqref{eq:driftless} is (globally) controllable.

By the formula of conservation of mass, \nota{this} transport problem under dynamical constraints is equivalent to the controllability of the homogeneous continuity equation
\begin{equation}\label{eq:continuity}
 \partial_t \mu_t + \mathrm{div}_g( \sum_{i=1}^{m} u_{i}f_{i}(x) \cdot\mu_t) = 0,
\end{equation}
where $\mathrm{div}_g$ denotes the divergence defined by the Riemannian volume form on $M$.
Namely, given $\mu$ and $\nu$ probability measures on $M$, find controls $u_1, \dots, u_m$ such that the solution, in the distributional sense, $\mu_t$ of~\eqref{eq:continuity} with initial condition $\mu_0 = \mu$ satisfies $\mu_1=\nu$.

The control of the  continuity equation boasts many possible applications ranging from crowd
dynamics~\cite{colombo2011control,CPT14} to opinion formation processes~\cite{bellomo2013}, from herd analysis~\cite{TB04swarming} to autonomous
multi-vehicle navigation~\cite{ren2008book} and, in general, it can be used to model uncertainties on the initial state of a deterministic controlled system~\cite{MQ18Mayer,CQ08Survey,JeanZidani22Ifac,JeanZidani22Hal}.
Indeed this controllability problem can be seen as a continuous analogous of  
``ensemble controllability''\nota{,} i.e. a problem involving a large number of more or less
identical subsystems that are being manipulated by a single source of command
signals~\cite{brockett2000stochastic,agrachev2016ensemble}.
\nota{The topic has gathered a fast growing interest in recent years due to the increased visibility
of the optimal control approach to multilayered artificial neural 
networks. 
Looking at the neural network as a  continuum of layers and each layer weight as a control function, 
Deep Learning can be seen as a simultaneous control problem for Neural ODEs or, equivalently, when the size of the network is large, as a control problem for Neural transport 
equations, see for instance~\cite{Andreis,scagliotti2023deep,rbzuazua2023,rbzuazua2024}.}

The idea is to control the evolution of a large number of non-interacting \emph{agents} or \emph{particles}. By adopting a \emph{macroscopic} viewpoint,
the collection of all agents can be represented as a mass evolving over time. If  the total number of agents is constant, then, up to a renormalization,  the total mass of the system is constantly equal to $1$. 
\nota{The evolution of this mass is described by a curve $t\mapsto \mu_t$, where $\mu_t$ is a probability measure representing the distribution of the system
at time $t$. If any single agent follows the trajectory given by \eqref{eq:driftless}, the curve $t\mapsto \mu_t$ is the solution of~\eqref{eq:continuity}.} 

In \nota{the first main result of} this paper we prove that, under suitable regularity conditions, the controllability of \eqref{eq:driftless}
by means of open-loop controls, is a sufficient condition for the controllability of~\eqref{eq:continuity} by means of time-varying feedback controls (Proposition~\ref{thm:controllability} below). \nota{Then, in } Theorem~\ref{maintheorem} below, we show that the optimal transport map, for the \nota{quadratic cost} (see Definition~\ref{def:wasserstein} below), between two sufficiently regular probability measures \nota{on the Euclidean space}, is the flow at time $1$ of system~\eqref{eq:driftless}. 

\nota{While it is more natural to study controllability on compact manifolds without boundary, the regularity of the optimal transport map heavily depends on the geometry of the ambient space, which can, in certain cases, obstruct its realization as the flow of a smooth dynamical system.
However, as it will be discussed in more detail in Section~\ref{sec:opt}, in the Euclidean setting there is a well-established literature on conditions for the optimal transport map to be a smooth diffeomorphism. On the other hand, the wild structure of the group of diffeomorphisms of a noncompact manifold or of a manifold with boundary leads to further challenges for the analysis of flows and diffeomorphisms.}

\subsection{Contents}

\nota{The paper is structured as follows. In Section~\ref{sec:cce}, we prove the finite-time controllability of the continuity equation on a compact Riemannian manifold $M$, under the sole assumption that the underlying dynamic is controllable. The proof assumes smooth regularity on the measures and relies on two auxiliary results: the first is a classical theorem by Moser, Theorem~\ref{moser}, which guarantees the existence of a diffeomorphism coupling the initial and the target measures; the second, Theorem~\ref{ac}, is a controllability result on the group of diffeomorphisms, which ensures the existence of time-dependent feedback controls that realize Moser's transport map as the time-one flow of the control system.

Section~\ref{sec:opt} is dedicated to the realization of the optimal transport map between two probability measures on a bounded domain $\Omega \subset \mathbb{R}^n$, as the time-one flow of a control system. These measures are assumed to be absolutely continuous with respect to the Lebesgue measure, with smooth and strictly positive densities. We begin with a brief overview of optimal transport preliminaries in Subsection \ref{prelopt}, and proceed to present the main result in Subsection \ref{mainresultsec:3}. In the proof, we first apply Lemma~\ref{prop:scag} to show that, under suitable conditions on $\Omega$, the optimal transport map $T$ exists and is an orientation-preserving diffeomorphism. Then, in the spirit of Theorem~\ref{ac}, we construct controls that realize $T$ in Proposition~\ref{prop:controllability}. The main technical challenge here stems from the lack of compactness; to address this, the strategy, developed in Subsection~\ref{proofprop}, is based on assuming that the control system is controllable on a larger open set $E$ containing $\Omega$.}

The rest of this section contains a review of the literature on the subject and the notation useful in the statements of the main results.

\subsection{State of the art}

The problem of the controllability of \eqref{eq:continuity} has been proposed by R.W. Brockett in~\cite{brockett2007optimal} as an ensemble controllability problem applied to robotics  (see also \cite{brockett2012notes}). Brockett studied the 
density equation associated with linear systems. 
A result, phrased in terms of controlling diffeomorphisms, is published in~\cite{agrachevcaponigro}.  The result is valid for driftless control-affine systems on compact manifolds, under the sole assumption of controllability of the system. \nota{In~\cite{Andreis},  the authors prove approximate controllability for the the class of diﬀeomorpisms isotopic to the identity using open-loop controls, under the additional hypothesis of \emph{Lie Algebra strong approximating property}}. This result has then been used in~\cite{scagliotti} to approximate the optimal transport map between two smooth densities with the flow of system~\eqref{eq:driftless}. 

The first result linking the controllability of the continuity equation and optimal transport is due to Khesin and Lee~\cite{khesinlee} (see also~\cite{agrachevlee}). They prove that, under the Lie Algebra rank condition, the optimal transport map is  the flow of a driftless control-affine system associated with time-varying feedback controls. Our result, Proposition~\ref{thm:controllability} below, in fact extends this result to controllable systems without any assumption on the Lie brackets.

Beside controllability, we also mention  optimal controllability results for~\eqref{eq:continuity} with open-loop controls:
in~\cite{borzi2019}
for a particular class of vector fields of the form
$
F= a(t, x) + u_1(t) + x u_2(t),
$
and in~\cite{pogodaev2016optimal} for optimal control without terminal constraints for nonlinear systems. Finally let us cite  results of controllability with localized vector fields, e.g.~\cite{DMR19,DMR20}.
We also refer to the recent preprint~\cite{rag2024} for a similar yet independent result.

\subsection{Preliminaries} 

Let $(M,g)$ be an $n$-dimensional oriented Riemannian manifold without boundary  and consider a driftless control affine system on $M$
\begin{equation}\label{sistemaprel}
\dot{x}(t) = \displaystyle \sum_{i=1}^{m} u_{i}f_{i}(x(t)), 
\end{equation}
where $\mathcal{F}=\{f_{1}, \dots , f_{m}\} \subset \mathrm{Vec}(M)$ is a family of smooth, i.e. $C^\infty$, vector fields on $M$, 
and $u=(u_1, \ldots,u_m)$ is a $m$-tuple of scalar controls.

If $M$ is not compact, we assume the vector fields $f_{1}, \dots , f_{m}$  to satisfy some growth condition in order to ensure global existence and uniqueness of solution of the Cauchy problem 
\begin{equation}\label{eq:cauchy}
\begin{cases}
\dot{x}(t) &= \displaystyle \sum_{i=1}^{m} u_{i}(t)f_{i}(x(t)),  \quad  \mbox{ a.e. } t\in [0,1], \\
x(0)&=x_{0},
\end{cases}
\end{equation} 
for any $x_{0} \in M$, in the Carath\'eodory sense (see e.g. \cite[Theorem 5.3]{hale}), associated with open-loop controls $u_1(t), \ldots,u_m(t) : [0,1] \to \R$ measurable bounded.
For instance, in the non-compact case, we assume that the vector fields $f_{1}, \dots , f_{m}$ satisfy $f_{j}(x) \leq \alpha(1+|x|)$, for some  $\alpha$ constant.

\begin{definition}
We say that \eqref{sistemaprel} is \emph{controllable} in $M$ if for any pair of points $x_{0}, x_{1}$ in $M$, there exists an open-loop control $u=(u_1,\dots,u_m)$, with $u_i:[0,1]\rightarrow \R$ measurable and bounded for any $i=1,\dots,m$, such that the associated solution $x(\cdot)$ of~\eqref{eq:cauchy}, with initial condition $x_0$, satisfies  $x(1) = x_1$.
\end{definition}

Throughout the paper, by a smooth map, we mean a $C^{\infty}$ map and by a diffeomorphism, we mean a smooth invertible map with smooth inverse. We denote by $\mathrm{Diff}(M)$, the infinite dimensional topological group of diffeomorphisms from $M$ to $M$ \nota{and} by $C^{\infty}(M,\R)$ the space of smooth functions from $M$ to $\R$. \nota{Both $\mathrm{Diff}(M)$ and $C^{\infty}(M,\R)$ are endowed with the topology of uniform convergence of the partial derivatives of any order on any compact subset of $M$, (see e.g. \cite{rudin}). Finally, by }$\Diff_{0}(M)$ \nota{we denote} the connected component of the identity in $\mathrm{Diff}(M)$. 
We observe that if $E$ is an open subset of $\R^n$, $\Diff_{0}(E)$ coincides with the group of orientation-preserving diffeomorphisms of $E$, whereas this is not true for $\Diff_{0}(M)$ where $M$ is a generic oriented smooth manifold, (see e.g. \cite{banyaga}).

Given a vector field $f\in\mathrm{Vec}(M)$, we denote by $e^{tf}$, the flow on $M$ at time $t$ generated by $f$, where $t \in \R$,  and $\{e^{tf}\mid t\in \R\}$ is a one parametric subgroup of $\Diff_{0}(M)$. Moreover, for any $P \in \mathrm{Diff}(M)$ and any $f\in\mathrm{Vec}(M)$, $P_*f$ denotes the pushforward of $f$ by $P$. Using the notation of chronological calculus (see e.g. \cite[Section 2.5]{agrachev2013control}), we also write $P_*f=P^{-1}\circ f \circ P.$
\begin{definition}\label{pcalli}
    Let $\mathcal{F}\subset \mathrm{Vec}(M)$ be the family of smooth vector fields given by system \eqref{sistemaprel}. We define $\mathcal{P}\subseteq\diff_0(M)$ as the set of all possible finite compositions of exponentials of vector fields in $\mathcal{F}$ rescaled by smooth functions, i.e. 
$$
\mathcal{P}:=\{e^{\alpha_kf_k}\circ\dots\circ e^{\alpha_1f_1}\mid k\in\N, \alpha_i \in C^{\infty}(M,\R), f_i \in \mathcal{F}\}.
$$
\end{definition}

In order to distinguish between open-loop controls and time-varying feedback controls, we denote the latter with bold $\mathbf{u} = (u_1(t,x), \ldots,u_m(t,x))$.

Given a Borel probability measure $\mu\in \mathscr{P}(M)$ and a time-varying feedback control $$\mathbf{u} = (u_1(t,x), \ldots,u_m(t,x)),$$ measurable bounded with respect to $t$ and smooth with respect to $x$, we consider the Cauchy problem in $\mathscr{P}(M)$
\begin{equation}\label{eq:cauchycontinuity}
\begin{cases} 
 \partial_t \mu_t & + \mathrm{div}_g\left(\sum_{i=1}^{m} u_{i}(t, \cdot)f_{i}\cdot\mu_t\right) = 0, \\
 \mu(0) &=\ \mu,
\end{cases}
\end{equation}
where the solution is meant in the distributional sense \nota{and $\mathrm{div}_g$ is the divergence induced by the Riemannian volume form $dV_g$ on $(M,g)$. In local coordinates $dV_g = \sqrt{|g|}\mathrm{d}x^1\wedge\dots\wedge \mathrm{d}x^n,$ where $\mathrm{d}x^i$ are $1$-forms that form a positively oriented basis for the cotangent bundle of the manifold.}

A sufficient condition for the existence and uniqueness of solution for the Cauchy problem~\eqref{eq:cauchycontinuity}, is the existence and uniqueness of solution for system~\eqref{eq:cauchy} (see, for instance,~\cite{AGS05} or ~\cite[Sections 4 and 5]{santambrogio2015optimal}).
Moreover, denoting by $\phi^t$ the flow at time $t$ of~\eqref{sistemaprel}, the solution of~\eqref{eq:cauchycontinuity} is given by 
\begin{equation}\label{eq:solution}
\mu_t = \phi^t_{\#} \mu,
\end{equation}
namely, for any Borel set $A \subset M$
$$
\phi^t_{\#} \mu(A) := \mu((\phi^{t})^{-1}(A)).
$$
When the initial datum $\mu$ of~\eqref{eq:cauchycontinuity} is absolutely continuous with respect to the Riemannian volume form  $dV_g$, then so is the 
 solution $\mu_t$ of~\eqref{eq:cauchycontinuity} for every $t\in[0,1]$. Namely, if 
$\mu=\rho_{0}dV_g$ for some smooth function $\rho_0: M \to \R$, then the solution is $\mu_t = \rho_t dV_g$, where $\rho_t: M \to \R$ is smooth for every $t$ and satisfies the \emph{density equation}
\begin{equation}\label{eq:density}
\rho_t(x) = \frac{1}{\mathrm{det}\nabla\phi^t(x)} \rho_0\left((\phi^{t})^{-1}(x)\right),
\end{equation}
for almost every $t\in[0,1]$ and for every $x\in M$.

\section{Controllability of the continuity equation}\label{sec:cce}

The main result of this section is a controllability result for the Continuity Equation in finite time, in the sense that every pair of measures with smooth  positive densities can be joined by an admissible solution of~\eqref{eq:cauchycontinuity} under the sole assumption of  controllability of the underlying dynamics on $M$.

We restrict our attention to the case in which $M$ is compact. Under this assumption, existence and uniqueness of solution of the Cauchy problem \eqref{eq:cauchy} \nota{are} ensured without any additional hypothesis on the vector fields $f_1,\dots,f_m$. Moreover, in the framework of absolutely continuous measures with smooth positive densities, controllability of the continuity equation \eqref{eq:cauchycontinuity} is equivalent to controllability of the density equation \eqref{eq:density}. 
\begin{definition}
    For any time-varying feedback control $\mathbf{u}=(u_1(t,x),\ldots,u_m(t,x))$,  measurable bounded with respect to $t \in [0,1]$ and smooth with respect to $x \in M$, we define $\phi_{\mathbf{u}}^t$ as the flow at time $t$ of \eqref{sistemaprel} associated with $\mathbf{u}$.  
\end{definition}

\begin{proposition}\label{thm:controllability} 
Let $(M,g)$ be a compact, $n$-dimensional, oriented Riemannian manifold without boundary \nota{and let $dV_g$ be the Riemannian volume form on $(M,g)$.} Let $\mu=\rho_{\mu}dV_g$ and $\nu=\rho_{\nu}dV_g$ be probability measures on $M$ with smooth and strictly positive densities $\rho_{\mu}$ and $\rho_\nu$. 

Then, if~\eqref{sistemaprel} is controllable in $M$, there exists a time-varying feedback control $\mathbf{u}$, measurable bounded with respect to $t \in [0,1]$ and smooth with respect to $x \in M$, such that $\nu$ is the $1$-time solution of the continuity equation \eqref{eq:cauchycontinuity}, associated to $\mathbf{u}$ with initial datum $\mu$.
\end{proposition} 

The main tools used \nota{to prove} Proposition~\ref{thm:controllability} are Moser's Theorem~\cite{moser,dacorognamoser} and \nota{a result on controllability on the group of diffeomorphisms in}~\cite{agrachevcaponigro}. Both results are stated below  for completeness.

\begin{theorem}\label{moser}[Moser~\cite{moser}]
Let $M$ be an $n$-dimensional, compact, smooth manifold without boundary. If $\sigma, \tau$ are two volume forms on $M$, then there exists $P \in \diff_0(M)$ such that $$\sigma=P_{\#}(\lambda \tau),$$
    where $$\lambda=\displaystyle \int_M \tau \Big/{\int_M \sigma}.$$
\end{theorem}

\begin{theorem}\label{ac}[Theorem 1.1~\cite{agrachevcaponigro}]
    Let $M$ be an $n$-dimensional, compact, smooth manifold without boundary. If \eqref{sistemaprel} is controllable in $M$, then there exists a neighborhood $\mathcal{O}$ of the identity in $\diff_0(M)$ and a positive integer $k$ such that every $P \in \mathcal{O}$ can be represented in the form $$P=e^{a_kf_k}\circ\dots\circ e^{a_1f_1},$$
    for some $f_1,\dots,f_k \in \mathcal{F}$ and $a_1,\dots,a_k \in C^{\infty}(M)$.  
\end{theorem}

\begin{remark} When studying lifting of trajectories on the group of diffeomorphisms, it is natural to
work with time-varying feedback controls. Indeed, if the $u_i$'s are
continuous feedback controls not depending on time, even the realization of 
locally asymptotically
stable equilibria for a control system $\dot x = f(x,u)$ is, in general, not possible~\cite{Brockett1983}. Hence the use of  time-varying
feedback controls has been introduced by J.-M. Coron (see Coron~\cite{Coron1,Coron2} or \cite[Section 11.2]{CoronBook}).
\end{remark}

\begin{remark}\label{remark_p}
    For any $P\in\mathcal{P}$ there exists a time-varying feedback control  $\mathbf{u}(t,x)=(u_1(t,x),\dots,u_m(t,x))$, smooth w.r.t. $x \in M$ and piece-wise constant w.r.t $t\in[0,1]$, such that $P=\phi^1_{\mathbf{u}}$, where $\phi^1_{\mathbf{u}}$ is the flow of system \eqref{sistemaprel} at time $1$ associated to $\mathbf{u}$. In fact, if $P \in \mathcal{P}$, we can write $P=e^{\alpha_kf_k}\circ\dots\circ e^{\alpha_1f_1}$ for suitable $k \in \N$, $\alpha_1,\dots,\alpha_k\, \in\, C^{\infty}(M,\R)$ and $f_1,\dots,f_k\,\in\,\mathcal{F}$. Thus it suffices to take the control $\mathbf{u}=(u_1,\dots,u_m)$, where for any $i=1,\dots,m$ 
    \begin{equation*}
        u_i(t,x):=\begin{cases}
            \alpha_j(x)  &\mbox{ if } f_j=f_i, \\
            0 &\mbox{ if } f_j\neq f_i,
        \end{cases}
    \end{equation*}
    for $x \in M$, $t \in [\frac{j-1}{k},\frac{j}{k})$ and $j=1,\dots,k$.
\end{remark}

\begin{proof}[Proof Proposition \ref{thm:controllability}]
Since $\rho_{\mu}$ and $\rho_{\nu}$ are smooth and strictly positive on $M$, $\mu$ and $\nu$ are volume forms on $M$. Moreover $\int_M \mu\,=\,\int_M \nu$ implies $\lambda=1$ in Theorem \ref{moser}. Hence it's possible to find $P\in \diff_0(M)$ such that $$\nu=P_{\#}\mu.$$
Thanks to Theorem \ref{ac}, the controllability assumption on system \eqref{sistemaprel} provides $\mathcal{O}\subseteq\diff_0(M)$ open neighborhood of the identity such that $\mathcal{O}\subseteq\mathcal{P}$.
\
Since $\diff_0(M)$ is a path-connected topological group, it is generated by any neighborhood of the identity. In particular, given $P\in\mathrm{Diff}_0(M)$, there exists a path $\{P^t\mid t\in[0,1]\}\subset\mathrm{Diff}_0(M)$ such that $P^0=\mathrm{Id}$ and $P^1=P$. For every $N \in \N$, consider the diffeomorphism 
$$
P^{k/N}\circ\big(P^{(k-1)/N}\big)^{-1}, \qquad k=1,\dots,N.
$$
If $N$ is large enough, $P^{k/N}\circ\big(P^{(k-1)/N}\big)^{-1} \in \mathcal{O}\subseteq \mathcal{P}$ for any $ k=1,\dots,N$, thus $$P=P\circ\big(P^{(N-1)/N}\big)^{-1}\circ P^{(N-1)/N}\circ\dots\circ P^{1/N}$$
and $P \in \mathcal{P}$. Then conclusion follows from Remark \ref{remark_p}, indeed $P=\phi^1_{\mathbf{u}}$ for some time-varying feedback control  $\mathbf{u}(t,x)=(u_1(t,x),\dots,u_m(t,x))$, smooth w.r.t. $x \in M$ and piece-wise constant w.r.t $t\in[0,1]$, and $$\nu=\phi^1_{{\mathbf{u}}_{\#}}\mu,$$
i.e. $\nu$ is the $1$-time solution of $\eqref{eq:cauchycontinuity}$ associated to $\mathbf{u}$ with initial datum $\mu$.
\end{proof}

\section{Optimal transport via time-varying feedback} \label{sec:opt}
This section is devoted to the problem of realizing the optimal transport map between two given probability measures as the flow of~\eqref{sistemaprel}. 
Since for any time-varying feedback control $\mathbf{u}(t,x)=(u_1(t,x),\dots,u_m(t,x))$, smooth w.r.t. $x$ and piece-wise constant w.r.t $t$, the flow of system \eqref{sistemaprel} associated to $\mathbf{u}$ is a diffeomorphism isotopic to the identity, we consider the particular case in which also the optimal transport map has this kind of regularity. However, in the Riemannian setting, even when considering measures that are absolutely continuous with smooth and strictly positive densities, the geometry of the manifold plays a crucial role on the regularity of the optimal transport map, see e.g. \cite{ma2005regularity}. As observed by Caffarelli \cite{caffarelli96}, also in the Euclidean case, one cannot expect any general regularity result without making some geometric assumptions on the support of the target measure. More precisely, \cite{caffarelli96} proved that if $\Omega$ is a bounded smooth open set in $\R^n$ and $\mu=\rho_{\mu}\mathscr{L}^n$, $\nu=\rho_{\nu}\mathscr{L}^n$ are two probability measures on $\Omega$ absolutely continuous w.r.t the Lebesgue measure, with $\rho_{\mu}$ and $\rho_{\nu}$ smooth and \nota{bounded away from zero}, if we assume that $\Omega$ is uniformly convex, then the optimal transport map $T:\Omega \rightarrow \Omega$ is a smooth diffeomorphism. Moreover ~\cite[Proposition 3.2]{scagliotti} ensures that $T$ is isotopic to the identity.
 
Therefore, in this setting, the possibility of realizing $T$ as the $1$-time flow of system \eqref{sistemaprel} follows from the controllability property of diffeomorphisms in $\Omega$ stated in Proposition \ref{prop:controllability} below. 
In Proposition \ref{prop:controllability} we need to adapt the result of \cite{agrachevcaponigro} to the non-compact case. To tackle this challenge, we assume that the closure of $\Omega$ is contained in a wider open set $E$ in which system $\eqref{sistemaprel}$ is controllable and we study the problem of realizing diffeomorphisms as flows of $\eqref{sistemaprel}$ restricted to $\Omega$.

We now introduce some basic notion of optimal transportation useful to state the main result. For further details we refer to the books~\cite{villani2009book},~\cite{santambrogio2015optimal}.
\subsection{Preliminaries on optimal transport}\label{prelopt}
\begin{definition}
Given a separable metric space $(X,d)$ and $\mu, \nu  \in \mathscr{P}(X)$, a Borel-measurable map 
$T: X \to X$ is called a \emph{transport
map} from $\mu$ to $\nu$ if $T_{\#} \mu = \nu$, namely,
$$
\nu(A) = \mu(T^{-1}(A)),
$$
for any Borel $A \subset X$.
The measure $T_{\#} \mu$ is called \emph{pushforward measure} of $\mu$ through the map $T$.
\end{definition}

The pushforward measure  satisfies
\begin{equation*}
\int_{X}\varphi(x)\spazio dT_{\#}\mu(x) \spazio \spazio = \spazio \spazio \int_{X}\varphi \circ T(x) \spazio d\mu(x),
\end{equation*}
for any test function $\varphi:X\to \R$, \nota{continuous} with compact support.\\
We denote by $\mathscr{P}_{2}(X)$ the set of Borel probability measures on $X$ having finite second moment, namely
\begin{equation*}
\mathscr{P}_{2}(X):=\Big\{\mu \in \mathscr{P}(X) \spazio : \spazio \int_X d(x,x_0)^2 d\mu(x) < +\infty \spazio < \spazio +\infty\Big\},
\end{equation*}

\nota{for any $x_0 \in X$.}

Finally, for any pair of probability measures $\mu, \nu \in \mathscr{P}(X)$, we define the set of \emph{admissible transport plans} between $\mu$ and $\nu$ as 
\begin{equation*}
\textrm{Adm}(\mu,\nu):=\{\gamma \in \mathscr{P}(X \times X):  (\pi_{1})_{\#}\gamma = \mu, (\pi_{2})_{\#}\gamma = \nu\},
\end{equation*}
where $\pi_{1}, \pi_{2} : X \times X \rightarrow X$ are the canonical projections on the first and the second component, respectively. 

\begin{definition} For every pair of probability measures $\mu, \nu \in  \mathscr{P}_{2}(X)$, the $2$-\emph{Wasserstein distance} $W_{2}(\mu, \nu)$ between $\mu$ and $\nu$ is defined as
\begin{equation}\label{wass}
W^2_{2}(\mu, \nu) := \mathrm{inf} \Big\{ \int_{X \times X} d(x,y)^{2} d\gamma (x,y) \mid \gamma  \in \mathrm{Adm}(\mu, \nu)\Big\}. 
\end{equation}
\end{definition}
We denote by $\mathrm {Opt}(\mu, \nu)$ the set of admissible transport plans realizing the infimum in~\eqref{wass}:
\begin{equation*}
\mathrm{Opt}(\mu, \nu):= \Big\{\gamma \in \textrm{Adm}(\mu,\nu) \mid  \int_{X \times X} d(x,y)^{2} \spazio d\gamma (x,y) = W_{2}^{2}(\mu, \nu) \Big\}. 
\end{equation*}
\begin{definition}\label{def:wasserstein}
Given a Borel map $T:X \rightarrow X$, we say that $T$ is an \emph{optimal transport map} between $\mu, \nu \spazio \in  \mathscr{P}_{2}(X)$ for the quadratic cost if $\gamma_{T}:=(\mathrm{Id},T)_{\#} \mu \in \mathrm{Opt}(\mu, \nu)$. 
\end{definition}

We are now ready to state the main result of this section. From now on $\mathscr{L}^n$ denotes the Lebesgue measure in $\R^n$.
\subsection{Main result}\label{mainresultsec:3}

\begin{theorem}\label{maintheorem}
Let $E$ be an open set of $\R^n$ and let \eqref{sistemaprel} be controllable in $E$. Let $\Omega$ be a smooth, bounded, connected, uniformly convex, open subset of $\R^n$ such that $\overline{\Omega}\subset E$. Let $\mu, \nu \in \mathscr{P}_{2}(\Omega)$  such that
\begin{itemize}
\item[$(i)$]$\mu=\rho_{\mu}\mathscr{L}^n$ and $\nu=\rho_{\nu}\mathscr{L}^n$;
\item[$(ii)$] $\rho_{\mu}, \rho_{\nu} \in C^{\infty}(\overline{\Omega}, \R)$;
\item[$(iii)$]$\rho_{\mu}, \rho_{\nu}$   are strictly positive on $\overline{\Omega}$;
\end{itemize}
and let $T:\Omega \rightarrow \Omega$ be the optimal transport map between $\mu$ and $\nu$. Then there exists a time-varying feedback control $\mathbf{u}$, measurable w.r.t. $t \in [0,1]$ and smooth w.r.t. $x \in E$, such that $$T=\phi^1_{\mathbf{u}}|_{\Omega},$$ where $\phi^1_{\mathbf{u}}|_{\Omega}$ is the restriction to $\Omega$ of the $1$-time flow of \eqref{sistemaprel} associated to $\mathbf{u}$.
\end{theorem}

Theorem \ref{maintheorem} is a consequence of the following Lemma, whose proof is given in \cite[Proposition 3.2]{scagliotti}, and Proposition \ref{prop:controllability}. Due to the technical nature of the proof of Proposition \ref{prop:controllability}, it will be deferred to section \ref{proofprop}.

\begin{lemma}\label{prop:scag}
Let $\Omega$ be a smooth, bounded, connected, uniformly convex, open subset of $\R^n$. Let $\mu, \nu \in \mathscr{P}_{2}(\Omega)$ satisfy the assumptions $(i), (ii), (iii)$ of Theorem \ref{maintheorem}, let $T:\Omega \rightarrow \Omega$ be the optimal transport map between $\mu$ and $\nu$, then $T\in \Diff_0(\Omega)$. 
\end{lemma}

\begin{proposition}\label{prop:controllability}
Let $E$ be an open set of $\R^n$. Let $\Omega$ be a smooth, bounded, connected, open subset of $\R^n$ such that $\overline{\Omega}\subset E$.
If \eqref{sistemaprel} is controllable in $E$, then  for any $P \in \mathrm{Diff}_{0}(\Omega)$ there exists a time-varying feedback control $\mathbf{u}$ measurable w.r.t. $t \in [0,1]$ and smooth w.r.t. $x \in E$, such that 
$$
P = \phi^{1}_{\mathbf{u}}|_\Omega,
$$
where $\phi^{1}_{\mathbf{u}} \in \diff_0(E)$ is the $1$-time flow of \eqref{sistemaprel} associated with $\mathbf{u}$ and $\phi^{1}_{\mathbf{u}}|_\Omega$ is its restriction to $\Omega$.
\end{proposition} 

\subsection{Proof of Proposition \ref{prop:controllability}}\label{proofprop}
In the proof of Proposition~\ref{prop:controllability} we retrace the steps of~\cite{agrachevcaponigro}. In order to deal with the lack of compactness we assume that $\Omega$ is a bounded domain of $\R^n$  compactly contained in a wider domain $E$, in which the system is controllable. 

The proof is based on the Orbit Theorem of Nagano-Sussmann~\cite{sussmann} and the Palis-Smale Lemma~\cite{PS70,PS00}. We restate both the results below for convenience to the reader.

\begin{definition}
    Given $P\in\mathrm{Diff}(\Omega)$, we define the support of $P$ as the closure of the set in which $P$ is different from the identity, that is 
    $$
    \mathrm{supp}\,P:=\overline{\{x\in \Omega\mid P(x)\neq x\}}.
    $$
\end{definition}

\begin{theorem}\label{orbittheorem}[Orbit Theorem~\cite{sussmann}]
If \eqref{sistemaprel} is controllable in $E$, then for any $x \in E$, there exist $Q_1,\dots,Q_n \in \{e^{t_lf_l}\circ \dots \circ e^{t_1f_1}\mid l\in\N, \, t_s\in\R, \, f_s \in \mathcal{F}, s=1,\dots,l\}$ and $f_1,\dots,f_n \, \in \, \mathcal{F}$ such that the vector fields 
\begin{equation}\label{campi}
    X_i:=Q_{i_*}f_i  \in \mathrm{Vec}(E), \quad i=1,\dots,n,
\end{equation} 
satisfy $$\mathrm{span}\{X_1(x),\dots,X_n(x)\}=\R^n.$$
\end{theorem}

\begin{lemma}\label{palissmale}[Palis-Smale~\cite{PS70}]
Let $\displaystyle\bigcup_{j=1}^{k} V_j= \Omega$ be a finite open covering of $\Omega$ ad let $\mathcal{O}$ be an open neighborhood of the identity in $\diff_0(\Omega)$. Then, for any $P \in \diff_0(\Omega)$, there exist $P_1,\dots,P_k$ diffeomorphisms in $\mathcal{O}$ such that $\mathrm{supp}\,P_j\subset V_j$ for any $j=1,\dots,k$ and $$P=P_k\circ\dots\circ P_1.$$
\end{lemma}
\nota{Theorem \ref{orbittheorem} and Lemma \ref{palissmale}} are well-known results and their proofs can be found in literature, e.g. \cite[Section 5.4]{Andreis} for the first one and \cite[Lemma 5.4]{agrachevcaponigro} for the second one. 
\begin{definition}
    Given $V$ open subset of $E$ and $p \in V$, we define $C_p^{\infty}(V,\R^n)$ as the set of smooth functions fixing $p$,  $$C_p^{\infty}(V,\R^n):=\{F \in C^{\infty}(V,\R^n)\mid F(p)=p\},$$
    endowed with the standard $C^{\infty}$ topology.
\end{definition}
\nota{

\begin{definition}\label{foliation}
    Given a vector field $W\in\mathrm{Vec}(\R^n)$, and an open set $V\subseteq \R^n$, a smooth map $\varphi:V\rightarrow \R^n$ preserves the $1$-foliation generated by the trajectories of the equation $\dot{\xi}=W(\xi)$ if, for any $q \in V$, there exists some $a=a(q)\in \R$ such that $\varphi(q)=e^{a(q)W}(q)$.
\end{definition}
}

Next result states that every diffeomorphism sufficiently close to the identity can be decomposed into $n$ diffeomorphisms, each preserving a direction. The result follows~\cite[Lemma~4.2]{agrachevcaponigro} in its proof; here, however, we present a slightly more general statement. In fact, the result applies to all compact sets of $\R^n$,  in which there is a linearly independent family of $n$ vector fields. 


\begin{lemma}\label{lemma:ift}
Let $V$ be a bounded open subset of $\R^n$ and let $X_1,\dots,X_n\,\in\,\mathrm{Vec}(\R^n)$ be linearly independent on $\overline{V}$, then for any $p\in V$ there exists $\mathcal{U}_p$ open neighborhood of the identity in $C_p^{\infty}(V,\R^n)$ such that for any $F \in \mathcal{U}_{p}$ exist $\varphi^{F}_1, \dots, \varphi^F_{n} \in C_p^{\infty}(V,\R^n)$ verifying
$$
F=\varphi^{F}_n\circ\dots\circ\varphi^{F}_1,
$$ 
where, for any $k=1,\dots,n$, $\varphi^{F}_k$ preserves the $1$-foliation generated by the trajectories of the equation $\dot{\xi}=X_k(\xi)$. 
\end{lemma} 

\nota
{
\begin{proof}
For any $q \in \overline{V}$ and for any ${t} = (t_1,\dots, t_n) \in \R^n$ consider  the map
    \begin{equation}\label{acca}
        H: (q, t) \mapsto e^{t_nX_n}\circ\dots\circ e^{t_1X_1}(q),
    \end{equation}
and call $H_q(\cdot):=H(q,\cdot)$. We denote by $DH_q(t)$ the differential of $H_q$ evaluated in $t$, hence 
$DH_q(t) = \partial_t H(q,t)$.
Since  $X_1,\dots,X_n\,\in\,\mathrm{Vec}(\R^n)$ are linearly independent on $\overline{V}$,
the matrix $DH_q(0)$ with columns $X_1(q),\dots,X_n(q)$ is invertible for any $q \in \overline{V}$.
Let $\lambda>0$ be the minimum over $q\in \overline{V}$ of the smallest singular value of each of the matrices  $DH_q(0)$, then
$$
\|(DH_q(0))^{-1}\|_2\leq \frac{1}{\lambda} \quad \mbox{ for every } q\in\overline{V}.
$$ 

Since the map $\tau\mapsto DH_q(\tau)$ is smooth, and $H(q,t)$ is smooth with respect both variables,  there exists $L>0$ such that 
$$
\|DH_q(\tau)-DH_q(0)\|_2\leq L |\tau|,
$$
for any $q \in \overline{V}$.

Now if $|\tau|<\frac{\lambda}{2L}$, then $DH_q(\tau)$ is invertible and 
\begin{equation}
    \|(DH_q(\tau))^{-1}\|_2 \leq \frac{2}{\lambda}.
\end{equation}
Indeed, denoting by $I$ the identity,
\begin{align*}
    \|I - (DH_q(0))^{-1}(DH_q(\tau)) \|_2 & \leq \|(DH_q(0))^{-1}\|_2 \|DH_q(\tau)-D_tH_q(0)\|_2 \\
    & \leq \frac{1}{\lambda} L |\tau| < \frac{1}{2},    
\end{align*}
and 
$$
\|(DH_q(\tau))^{-1}\|_2 \leq \frac{\|(DH_q(0))^{-1}\|_2}{1-\|I - (DH_q(0))^{-1}(DH_q(\tau))\|_2} \leq \frac{2}{\lambda}.
$$

For the Inverse Function Theorem (see, for instance,~\cite[Appendix A]{quantIFT} for a quantitative statement) $H_q$ is a smooth diffeomorphism from $B_{\frac{\lambda}{2L}}(0)$ to $H_q(B_{\frac{\lambda}{2L}}(0))$ and, moreover, for any $\gamma \leq \lambda^2/4L$, 
$B_{\gamma}(q) \subseteq  H_q(B_{\frac{\lambda}{2L}}(0))$.
For any fixed $p \in V$ and $\gamma \leq \lambda^2/4L$, define
\begin{equation*}
   \mathcal{U}^{\gamma}_{p}  :=  \{F \in C_p^{\infty}(V,\R^n) \mid \Vert F - I\Vert_{C^1(V)} < \gamma\}.
\end{equation*}
By definition of $\mathcal{U}^{\gamma}_{p}$, for every $F \in \mathcal{U}^{\gamma}_{p}$ and every $q \in V$, $F(q)\in B_{\gamma}(q)$. 
Hence for every $F \in \mathcal{U}^{\gamma}_{p}$ and every $q \in V$ there exists a unique $n$-uple of real numbers $(t_1^F(q),\dots,t_n^F(q))=t^{F}(q)\in H_q^{-1}(B_{\gamma}(q))$ such that 
\begin{equation*}
    F(q)=e^{t_n^F(q)X_n}\circ \dots \circ e^{t_1^F(q)X_1}(q)=H_{q}(t^F(q)), 
\end{equation*}
where $t^F(q)$ is $C^\infty$ w.r.t. $q$ for the Implicit Function Theorem. Moreover  $(t_1^F(p),\dots,t_n^F(p))=(0,\dots,0)$  for any $F \in \mathcal{U}^{\gamma}_{p}$. From now on we fix $F\in \mathcal{U}^{\gamma}_p$ and, for the sake of readability,  we omit the superscript in $t(q):=t^F(q)$.

Consider, for $k=1,\dots,n$,   
$$
\Phi^F_k(q):=e^{t_k(q)X_k}\circ \dots \circ e^{t_1(q)X_1}(q).
$$
Let $\Phi^F_0=I$ and note that $\Phi^F_n=F$. For every $k$, $\Phi^F_k$ is a smooth map from $V$ to $\R^n$ and its differential at $q \in V$ is
\begin{align}
    D\Phi^F_k(q)(\xi)= & D(e^{\tau_kX_k}\circ \dots \circ e^{\tau_1X_1})(q)\xi + \nonumber \\
    & + \sum_{i=1}^{k}D\left(e^{\tau_kX_k}\circ\dots\circ e^{\tau_{i+1}X_{i+1}}\right)(q) \nabla t_i(q)\cdot\xi X_i (e^{\tau_{i}X_{i}}\circ\dots\circ e^{\tau_1X_1}(q)), \label{differenziale}
\end{align}
at $\tau_i = t_i(q)$ for $i=1,\ldots,k$ and for $\xi \in \R^n$.

Let us prove that $\|t\|_{C^1(V)}\rightarrow 0$ as $\gamma \rightarrow 0$.
Indeed  $\lim_{\gamma \to 0}\|t\|_{C^0(V)} = 0$ since $t(q)=H_q^{-1}(F(q))$ and
\begin{align*}
|t(q)| & = 
|H_q^{-1}(F(q))| \\
&=   |H_q^{-1}(F(q)) - H_q^{-1}(q)| \\ 
&\leq \sup_{|\tau| <\lambda/2L}\|(DH_q(\tau))^{-1}\|_2 
\|F - I\|_{C^0(V)} \\
&< \frac{2}{\lambda} \gamma.
\end{align*}

In order to estimate $\|Dt(q)\|_2$ we differentiate both sides of $F(q)=H(q,t(q))$ w.r.t. $q$ and obtain
$$
DF(q) = \partial_q H(q,t(q)) + \partial_t H(q,t(q)) \cdot Dt(q).
$$
Since  $\partial_t H(q,t(q))=DH_q(t(q))$ is invertible 
\begin{align*}
     \|Dt(q)\|_2&\leq \|(DH_q(t(q)))^{-1}\|_2 \Big(\|DF(q)-I\|_2+\|I-\partial_q H(q,t(q))\|_2\Big)\\
     &\leq \frac{2}{\lambda}\Big(\gamma+\|I-\partial_q H(q,t(q))\|_2\Big).
\end{align*}
Let $K\subset \R^n$ be a compact set containing all the flow lines of $e^{t_nX_n}\circ\dots\circ e^{t_1X_1}(q)$ with $q \in \overline{V}$ and $(t_1,\dots,t_n) \in B_{\frac{\lambda}{2L}}(0)$, let $M>0$ such that for any $k=1,\dots,n$ 
$\max_{\xi \in K}\|JX_i(\xi)\|<M,$
where $JX_i(\xi)$ is the Jacobian matrix of the vector field $X_i$ at $\xi$. 
By chain rule 
\begin{align*}
    \partial_q H(q,t(q))&=D\big(e^{\tau_nX_n}\circ\dots\circ e^{\tau_1X_1}\big)(q)\\
    &=De^{\tau_nX_n}(\Phi_{n-1}^F(q))\cdot De^{\tau_{n-1}X_{n-1}}(\Phi_{n-2}^F(q))\cdot \ldots \cdot De^{\tau_{1}X_{1}}(q), 
\end{align*}
 with $(\tau_1,\dots, \tau_n) = (t_1(q), \dots, t_n(q))$. 
Since for any $k=1,\dots,n$, the matrix $De^{\tau_kX_k}(\Phi_{k-1}^F(q))$ is the solution at time $\tau_k=t_k(q)$ of the differential equation 
    $$
    \frac{d}{ds}\eta=JX_k(\xi(s))\eta, \quad \eta(0)=I,
    $$
where $\xi(s)=e^{sX_k}(\Phi_{k-1}^F(q))$, we have
$$
De^{\tau_kX_k}(\Phi_{k-1}^F(q))=I+\tau_kJX_k(\Phi_{k-1}^F(q))+o(|\tau_k|),
$$
hence
$$
\|De^{\tau_kX_k}(\Phi_{k-1}^F(q)) -I\|\leq M |\tau_k|+o(|\tau_k|). 
$$
Thus, for the entire product $\|I-\partial_q H(q,t(q))\|_2\leq C|t(q)|\leq \frac{2C}{\lambda} \gamma$ , for some constant $C$ depending only on $M$ and $n$. This implies that
$
\sup_{q\in V}\|Dt(q)\|_2\leq \frac{2}{\lambda}(1+C) \gamma.
$
In particular $\|t^F\|_{C^1(V)}\rightarrow 0$ as $\gamma \rightarrow 0$,   for any $F \in \mathcal{U}^{\gamma}_p$.

By $\|t^F\|_{C^1(V)}\rightarrow 0$ as $\gamma \rightarrow 0$ and~\eqref{differenziale}, we can choose $\gamma:=\overline{\gamma}$ small enough to ensure $\|D_q\Phi^F_k-I\| <1$ for any $q \in V$ and any $F \in \mathcal{U}^{\overline{\gamma}}_p$. Therefore, denoting $\mathcal{U}_p:=\mathcal{U}^{\overline{\gamma}}_p$, we have that for any  $q \in V$ and any $F \in \mathcal{U}_p$, $D_q\Phi^F_k$ is invertible and by the Implicit Function Theorem, so is $\Phi^F_k:V\rightarrow \Phi^F_k(V).$ For any $F \in \mathcal{U}_p$, let us fix 
$
\varphi^F_k: q\mapsto e^{t_k^F((\Phi^{F}_{k-1})^{-1}(q))X_k}(q),
$
so that 
$$
F(q)=\varphi^{F}_n\circ\dots\circ\varphi^{F}_1(q)
$$
for every $q\in V$ and $\varphi_k^F$ preserves the $1$-foliation generated by the trajectories of the equation $\dot{q}=X_k(q)$.
\end{proof}
}

\begin{lemma}\label{lem:cuore} 
Let $V$ be a bounded open subset of $\R^n$ and let $X_1,\dots,X_n \in\mathrm{Vec}(\R^n)$ be linearly independent on $\overline{V}$. Then for any $p\in V$ there exists $\mathcal{U}_p$ open neighborhood of the identity in $C_p^{\infty}(V,\R^n)$ such that any $F \in \mathcal{U}_{p}$ can be written as 
$$
F=e^{a_nX_n}\circ\dots\circ e^{a_1X_1},
$$  
for $a_1, \dots, a_{n} \in C^{\infty}(V,\R)$ satisfying $a_1(p)=\dots=a_n(p)=0$.
\end{lemma}

\nota{
This result is a restatement of~\cite[Proposition 4.1]{agrachevcaponigro}. 
Since the proof is exactly as in~\cite[Proposition 4.1]{agrachevcaponigro} we provide only a sketch of the main steps of the proof. The only difference is in the application of Lemma~\ref{lemma:ift}. Indeed, as in Lemma~\ref{lemma:ift} the result actually yields in all neighborhood of $\R^n$ in which there is a linearly independent family of $n$ vector fields.
} The idea behind the proof is the following. Consider a linear ODE on $U \subset \R$ of the form 
\begin{equation*}
    \begin{cases}
        \dot{x}&=\beta x\\
        x(0)&=x_0. 
    \end{cases}
\end{equation*}
Since the solution is $x(t)=e^{t\beta}x_0$, the $1$-time flow of this equation is given by the linear diffeomorphism $\phi^1(x)=e^{\beta}x.$ Conversely any linear diffeomorphism of $U \subset \R$, $\phi(x)=\alpha x$, with $\alpha\neq 1$ and $\alpha>0$, is the exponential of the linear vector field $\mathrm{log}(\alpha)x\frac{\partial}{\partial x}$. As the problem is extremely simple in the one-dimensional case with linear flow, the strategy of the proof is to reduce it to this case. 

\nota{
\begin{proof}[Sketch of the Proof of Lemma~\ref{lem:cuore}]
Thanks to Lemma \ref{lemma:ift}  for any $p \in V,$ there exists $\mathcal{U}_p$ open neighborhood of the identity in $C_p^{\infty}(V,\R^n)$ such that for any $F \in \mathcal{U}_{p}$ exist $\varphi_1, \dots, \varphi_{n} \in C_p^{\infty}(V,\R^n)$ verifying
$$
F(q)=\varphi_n\circ\dots\circ\varphi_1(q)
$$ 
for any $q\in V$,
where, for any $k=1,\dots,n$, $\varphi_k$ preserves the $1$-foliation generated by the trajectories of the equation $\dot{\xi}=X_k(\xi)$. 

The problem then reduces to realize each $\varphi_k$ as the flow of $X_k$ multiplied by a suitable smooth function. Since for any $q \in V$, there exists $f(q) \in \R$ such that  $\varphi_k(q) = e^{f(q)X_k}(q)$ (see Definition~\ref{foliation} above), we can treat the problem as a one-dimensional problem along the direction $X_k$ with the remaining $n-1$ coordinates playing the role of parameters. For simplicity let us omit the subscript $k$ and the $n-1$ parameters. Moreover, up to a change of coordinates, let us assume that $X_k=\frac{\partial}{\partial x}$ and that $p=0$. For further details we refer to~\cite[Proposition 4.1]{agrachevcaponigro}.

\begin{itemize}
\item Given a diffeomorphism $\varphi$ in the connected component of the identity, let 
$\alpha = \log(\varphi'(0))$ and consider the homotopy
$$
\varphi_t(x) =
\begin{cases}
e^{\alpha(t-1)}\frac{\varphi(tx)}{t} &\quad t\in (0,1]\\
x &\quad t=0.
\end{cases}
$$

\item Let $a(t,x)\frac{\partial}{\partial x}$ be the vector field having $\varphi_t$ as flow, or, using the \emph{chronological notation} (see~\cite{agrachev2013control}),
$$
\varphi_t = \rexp{t}{a(\tau,\cdot) \frac{\partial}{\partial x}}{\tau}.
$$
Note that 
$
a(t,x) = \alpha x + b(t,x)x,$ with  $b(t,0)=0.$

\item  Consider the solution $u(t,x)$ of the first order linear PDE
$$
a(t,x)  \frac{\partial u}{\partial x} (t,x) + \frac{\partial u}{\partial t} (t,x) + b(t,x)u(t,x)=0.
$$

\item Consider time-depending change of coordinates $\psi_t(x) = x u(t,x)$, that linearizes the flow of $a(t,x)  \frac{\partial}{\partial x}$, namely
$$
\psi_t \circ \varphi_t \circ \psi^{-1}_t = e^{t\alpha x \frac{\partial}{\partial x}}.
$$
\end{itemize}
In conclusion, $\varphi$ is the $1$-time flow of the vector field
$
\left(\psi_1\right)_* \alpha x \frac{\partial}{\partial x}.
$
\end{proof}

\begin{lemma}\label{nostro}
        Let $\Omega$ open subset of $\R^n$ and $V \subset \Omega$, open and bounded. Let $X_1,\dots,X_n\,\in\,\mathrm{Vec}(\R^n)$ linearly independent on $\overline{V}$. Then there exists $\mathcal{U}$, open neighborhood of the identity in $\diff_0(\Omega)$, such that for any $P \in \mathcal{U}$ having $\mathrm{supp}\,P\subset V$ 
        $$P|_V= e^{s_nbX_n}\circ\dots\circ e^{s_1bX_1}\circ e^{a_nX_n}\circ\dots\circ e^{a_1X_1},$$
        for  $a_1,\dots,a_n \in C^{\infty}(V,\R)$, $s_1,\dots,s_n \in \R$ and $b \in C^{\infty}(\Omega,\R)$.
\end{lemma}
}

\begin{proof}
Consider a cut-off function $b \in C^{\infty}(\Omega;\R)$ such that $\mathrm{supp}\,b\subset V$ and such that $b^{-1}(1)$ has non-empty interior. Given $\varepsilon>0$ and $B_{\varepsilon}(0)$ in $\R^n$, we define the map 
$$
\Lambda:B_{\varepsilon}(0) \rightarrow \mathrm{Diff}_0(\Omega),
$$ 
associating with any $(s_1,\dots,s_n) \in B_{\varepsilon}(0)$ the diffeomorphism $\Lambda(s_1,\dots,s_n)$ defined by
$$
\xi \mapsto \Lambda(s_1,\dots,s_n)(\xi):=e^{s_nbX_n}\circ\dots\circ e^{s_1bX_1}(\xi). 
$$
Note that,  since $b(\xi)=0$ for any $\xi \in \Omega\setminus\overline{V}$, we have $\mathrm{supp}\,\Lambda(s_1,\dots,s_n)\subset V$. Let $p\in V$ in the interior of $b^{-1}(1)$ and $\varepsilon>0$ sufficiently small such that the map 
$$
\Lambda^{p}:(s_1,\dots,s_n)\mapsto \Lambda(s_1,\dots,s_n)(p),
$$ 
is a local diffeomorphism on $B_{\varepsilon}(0)$ mapping $(0,\dots,0)$ to $p$.
Thus $\Lambda^p(B_{\varepsilon}(0))$ is an open neighborhood of $p$ in $V$.

Let $\mathcal{U}$ be a neighborhood of the identity in $\Diff_0(\Omega)$ sufficiently small such that $P(p) \in \Lambda^p(B_{\varepsilon}(0))$ for any $P\in\mathcal{U}$ having $\mathrm{supp}\,P \subset V$ and consider $\mathcal{U}_p$ the open neighborhood of the identity in $C_p^{\infty}(V,\R^n)$ given by Lemma \ref{lem:cuore}. 

If  $\mathcal{U}$ and $\varepsilon > 0$ are small enough, then for any $\Phi \in \Lambda(B_{\varepsilon}(0))$ and for any $P \in\mathcal{U}$ with $\mathrm{supp}\,P \subset V$ such
that $\Phi\circ P(p)=p$, we have that $\Phi\circ P|_V \in \mathcal{U}_p$.

Now, given $P \in \mathcal{U}$ with $\mathrm{supp}\,P \subset V$, since $P(p) \in \Lambda^p(B_{\varepsilon}(0))$, there exists $(s_1,\dots,s_n)\in B_{\varepsilon}(0)$ such that 
$$ 
P(p)=e^{s_nbX_n}\circ\dots\circ e^{s_1bX_1}(p),
$$
thus 
$$ 
e^{-s_1bX_1}\circ\dots\circ e^{-s_nbX_n}\circ P(p)=p. 
$$
and $e^{-s_1bX_1}\circ\dots\circ e^{-s_nbX_n}\circ P |_V \in \mathcal{U}_p$. Therefore, by Lemma \ref{lem:cuore}, there exist $a_1,\dots,a_n \in C^\infty(V,\R)$ such that 
$$
e^{-s_1bX_1}\circ\dots\circ e^{-s_nbX_n}\circ P|_V =e^{a_1 X_1} \circ \dots \circ  e^{a_n X_n},
$$
that is 
$$
P|_V=e^{s_1 b X_1} \circ \dots \circ  e^{s_n b X_n} \circ  e^{a_1 X_1} \circ \dots \circ  e^{a_n X_n}.
$$
\end{proof}

\begin{proof}[Proof Proposition~\ref{prop:controllability}]
Let $\mathcal{F}=\{f_1,\dots,f_m\}\subset\mathrm{Vec}(E)$ and $\mathcal{P}:=\{e^{\alpha_kf_k}\circ\dots\circ e^{\alpha_1f_1}\mid k\in\N,\, \alpha_i \in C^{\infty}(M,\R),\, f_i \in \mathcal{F}\}$, as in Definition \ref{pcalli}. From Remark \ref{remark_p}, Proposition \ref{prop:controllability} follows from the inclusion 
\begin{equation}\label{inclusion}
    \diff_0(\Omega)\subseteq\mathcal{P}|_{\Omega},
\end{equation}
where $\mathcal{P}|_{\Omega} = \{P|_\Omega \, : \, P \in \mathcal{P}\}$.

Since the system is controllable in $E$, applying Theorem \ref{orbittheorem} to the points of $\overline{\Omega}\subset E$, we have that for any $x \in \overline{\Omega}$, there are $Q^x_1,\dots,Q^x_n \in \{e^{t_lf_l}\circ \dots \circ e^{t_1f_1}\mid l\in\N, \, t_s\in\R, \, f_s \in \mathcal{F}, s=1,\dots,l\}$ and $f^x_1,\dots,f^x_n \, \in \, \mathcal{F}$ such that the vector fields $X_i^x:=Q^x_{i_*}f^x_i$, for $i=1,\dots,n$, are linearly independent at $x$. 

Let $U_x\subset \subset E$ be a bounded open neighborhood of $x$, such that the vector fields $X_1^x,\dots,X_n^x$ are linearly independent on $\overline{U_x}$.

Since $\Omega$ is 
compactly contained in $E$, given the open covering $$\overline{\Omega}\subset \displaystyle\bigcup_{x \in \overline{\Omega}}U_x,$$ we can extract a finite open subcovering $$\overline{\Omega}\subset \displaystyle\bigcup_{j=1}^kV_j,$$ where $U_j \in \{U_x\mid x\in \overline{\Omega}\}$ for any $j=1,\dots,k$.

Now define $V_j:=U_j\cap \Omega$, thus 
\begin{equation}\label{finite:covering}
    \,\Omega\,=\displaystyle\bigcup_{j=1}^kV_j\,.
\end{equation}
By construction, for any of the open sets $V_j$, there is $x_j\in \overline{\Omega}$ such that the vector fields $X_1^{x_j},\dots,X_n^{x_j}$ are linearly independent over $\overline{V_j}$. For simplicity, let's rename these vector fields $X^j_1,\dots,X^j_n$.\\
For any of the open sets $V_j$, we consider $\mathcal{U}_j$ open neighborhood of the identity in $\diff_0(\Omega)$ given by Lemma \ref{nostro}, where we take $X_1^j,\dots,X_n^j$, as the family of linearly independent vector fields over $\overline{V_j}$.

Applying Lemma \ref{palissmale} with respect to the finite open covering of $\Omega$ given by \eqref{finite:covering} and with respect to $$\mathcal{O}=\displaystyle\bigcap_{j=1}^{k}\mathcal{U}_j,$$ we get that for any $P \in \diff_0(\Omega)$ there exist $$P_1,\dots,P_k\,\in\,\displaystyle\bigcap_{j=1}^{k}\mathcal{U}_j,$$ such that $\mathrm{supp}P_j\subset V_j$ and $$P=P_k\circ\dots\circ P_1.$$ 
In particular $P_j \in \mathcal{U}_j$ for any $j=1,\dots,k$. 
Thus, from Lemma~\ref{nostro}, there are $a_1^j,\dots, a_n^j \in C^{\infty}(V_j,\R)$, $b^j \in C^{\infty}(\Omega,\R)$ and $s_1^j,\dots,s_n^j \in \R$ such that $$P_j|_{V_j}=e^{s_n^jb^jX_n^j}\circ\dots\circ e^{s_1^jb^jX_1^j}\circ e^{a_n^jX_n^j}\circ\dots\circ e^{a_1^jX_1^j}.$$
Since $\mathrm{supp}\,P_j$ and $\mathrm{supp}\,b$ are compactly contained in the open set $V_j$, over which the vector fields $X_1^j,\dots,X_n^j$ are linearly independent, also $\mathrm{supp}\,a_i^j$, for $i=1,\dots,n$, are compact subsets of $V_j$. 
This allows us to define the smooth extensions $\overline{a}_i^j\,\in\,C^{\infty}(E,\R)$, where 
\begin{equation*}
    \overline{a}_i^j(x)= \begin{cases}
        a_i^j(x) & x\in V_j, \\
        0  &x\notin V_j,
    \end{cases}
\end{equation*}
and to extend $P_j$ to $\overline{P}_j\in \diff_0(E)$, where $$\overline{P}_j=e^{s_n^jb^jX_n^j}\circ\dots\circ e^{s_1^jb^jX_1^j}\circ e^{\overline{a}_n^jX_n^j}\circ\dots\circ e^{\overline{a}_1^jX_1^j}.$$
Thus $$P=\overline{P}_k\circ \dots \circ \overline{P}_1|_{\Omega}.$$
The conclusion follows from the definition of action of diffeomorphisms on vector fields, see e.g. \cite[Section 2.5]{Andreis}. Indeed for any $i=1,\dots,n$ and $j=1,\dots,k$, the corresponding $X_i^j$ is of the form $X_i^j=Q^j_{i_*}f^j_i$ for some $Q^j_i\in \{e^{t_lf_l}\circ \dots \circ e^{t_1f_1}\mid l\in\N, \, t_s\in\R, \, f_s \in \mathcal{F}, s=1,\dots,l\}\subset \mathcal{P}$ and $f^j_i \in \mathcal{F}$, thus 
$$ 
e^{s_ib^jX_i^j}=(Q_i^j)^{-1}\circ e^{(s_i\overline{b}^j\circ Q_i^j)f_i^j}\circ Q_i^j,
$$
and 
$$
e^{\overline{a}_i^jX_i^j}=(Q_i^j)^{-1}\circ e^{(\overline{a}_i^j\circ Q_i^j)f_i^j}\circ Q_i^j,
$$
implying $e^{s_ib^jX_i^j}\in \mathcal{P}$ and
$e^{\overline{a}_i^jX_i^j} \in \mathcal{P}$, hence 
$\overline{P}_k\circ \dots \circ \overline{P}_1\in \mathcal{P}$ and inclusion \eqref{inclusion} is proved. 
\end{proof}

\section*{ACKNOWLEDGEMENTS}                               
The authors are grateful to Alfio Borz\`i and Alessandro Scagliotti for the precious  discussions that inspired this work.

\section*{Declarations}
\noindent{\bf Ethical Approval} Not applicable.

\noindent{\bf Funding} This work is supported by the MIUR Excellence Department Project MatMod@TOV awarded to the Department of Mathematics, University of Rome Tor Vergata (CUP E83C23000330006), by 
the projects PRIN 2022-E53D23005910006 and PRIN 2022 PNRR- 
E53D23017910001 founded by the EU - Next Generation EU,
and by Gnampa - INdAM.

\bibliographystyle{plain}
\bibliography{bibliografiaopt2}

\end{document}